\documentclass[10pt]{amsart}

\usepackage[nobysame,abbrev,alphabetic]{amsrefs}
\usepackage{amssymb}
\usepackage[all]{xy}

\DeclareMathOperator{\cd}{cd}
\DeclareMathOperator{\Char}{char}
\DeclareMathOperator{\Gal}{Gal}

\DeclareMathOperator{\Img}{Im}

\DeclareMathOperator{\Ker}{Ker}

\DeclareMathOperator{\res}{res}

\DeclareMathOperator{\Sym}{Sym}
\DeclareMathOperator{\trg}{trg}

\begin{document}

\newtheorem{thm}{Theorem}[section]
\newtheorem{cor}[thm]{Corollary}
\newtheorem{lem}[thm]{Lemma}
\newtheorem{prop}[thm]{Proposition}
\newtheorem{defin}[thm]{Definition}
\newtheorem{exam}[thm]{Example}
\newtheorem{examples}[thm]{Examples}
\newtheorem{rem}[thm]{Remark}
\newtheorem{case}{\sl Case}
\newtheorem{claim}{Claim}
\newtheorem{prt}{Part}
\newtheorem*{mainthm}{Main Theorem}
\newtheorem*{thmA}{Theorem A}
\newtheorem*{thmB}{Theorem B}
\newtheorem*{thmC}{Theorem C}
\newtheorem*{thmD}{Theorem D}
\newtheorem*{thmE}{Theorem E}
\newtheorem{question}[thm]{Question}
\newtheorem*{notation}{Notation}
\swapnumbers
\newtheorem{rems}[thm]{Remarks}
\newtheorem*{acknowledgment}{Acknowledgment}

\newtheorem{questions}[thm]{Questions}
\numberwithin{equation}{section}

\newcommand{\ab}{\mathrm{ab}}
\newcommand{\Coker}{\mathrm{Coker}}
\newcommand{\dec}{\mathrm{dec}}
\newcommand{\dirlim}{\varinjlim}
\newcommand{\discup}{\ \ensuremath{\mathaccent\cdot\cup}}
\newcommand{\gr}{\mathfrak{gr}}
\newcommand{\nek}{,\ldots,}
\newcommand{\inv}{^{-1}}
\newcommand{\isom}{\cong}
\newcommand{\sep}{\mathrm{sep}}
\newcommand{\sym}{\mathrm{sym}}
\newcommand{\tagg}{^{''}}
\newcommand{\tensor}{\otimes}

\newcommand{\alp}{\alpha}
\newcommand{\gam}{\gamma}
\newcommand{\Gam}{\Gamma}
\newcommand{\del}{\delta}
\newcommand{\Del}{\Delta}
\newcommand{\eps}{\epsilon}
\newcommand{\lam}{\lambda}
\newcommand{\Lam}{\Lambda}
\newcommand{\sig}{\sigma}
\newcommand{\Sig}{\Sigma}
\newcommand{\bfA}{\mathbf{A}}
\newcommand{\bfB}{\mathbf{B}}
\newcommand{\bfC}{\mathbf{C}}
\newcommand{\bfF}{\mathbf{F}}
\newcommand{\bfP}{\mathbf{P}}
\newcommand{\bfQ}{\mathbf{Q}}
\newcommand{\bfR}{\mathbf{R}}
\newcommand{\bfS}{\mathbfS}
\newcommand{\bfT}{\mathbf{T}}
\newcommand{\bfZ}{\mathbf{Z}}
\newcommand{\dbA}{\mathbb{A}}
\newcommand{\dbC}{\mathbb{C}}
\newcommand{\dbF}{\mathbb{F}}
\newcommand{\dbN}{\mathbb{N}}
\newcommand{\dbQ}{\mathbb{Q}}
\newcommand{\dbR}{\mathbb{R}}
\newcommand{\dbZ}{\mathbb{Z}}
\newcommand{\grf}{\mathfrak{f}}
\newcommand{\gra}{\mathfrak{a}}
\newcommand{\grm}{\mathfrak{m}}
\newcommand{\grp}{\mathfrak{p}}
\newcommand{\grq}{\mathfrak{q}}
\newcommand{\calA}{\mathcal{A}}
\newcommand{\calB}{\mathcal{B}}
\newcommand{\calC}{\mathcal{C}}
\newcommand{\calE}{\mathcal{E}}
\newcommand{\calG}{\mathcal{G}}
\newcommand{\calH}{\mathcal{H}}
\newcommand{\calK}{\mathcal{K}}
\newcommand{\calL}{\mathcal{L}}
\newcommand{\calW}{\mathcal{W}}
\newcommand{\calV}{\mathcal{V}}

\title[Quotients of absolute Galois groups]{Quotients of absolute Galois groups which determine
the entire Galois cohomology}

\author{Sunil K.\ Chebolu}
\address{Department of Mathematics\\
Illinois State University\\
Campus box 4520\\
Normal, IL 61790\\
USA} \email{schebol@ilstu.edu}

\author{Ido Efrat}
\address{Mathematics Department\\
Ben-Gurion University of the Negev\\
P.O.\ Box 653, Be'er-Sheva 84105\\
Israel} \email{efrat@math.bgu.ac.il}

\author{J\'an Min\'a\v c}
\address{Mathematics Department\\
University of Western Ontario\\
London\\
Ontario\\
Canada N6A 5B7} \email{minac@uwo.ca}
\thanks{
Sunil Chebolu was supported by a New Faculty Initiative Grant at Illinois State
University.
Ido Efrat was supported by the Israel Science Foundation (grant No.\ 23/09).
J\'an Min\'a\v c was supported in part by  National Sciences and Engineering Council of Canada grant R0370A01.}

\keywords{absolute Galois group, Galois cohomology, descending central sequence, $W$-group}
\subjclass[2000]{Primary 12G05; Secondary 12F10, 12E30}

\begin{abstract}
For a prime power $q=p^d$ and a field $F$ containing a root of unity of order $q$ we show that the Galois cohomology ring
$H^*(G_F,\dbZ/q)$ is determined by a quotient $G_F^{[3]}$ of the absolute Galois group $G_F$ related to its
descending $q$-central sequence.
Conversely, we show that $G_F^{[3]}$ is determined by the lower cohomology of $G_F$.
This is used to give new examples of pro-$p$ groups which do not occur as absolute Galois groups of fields.
\end{abstract}

\maketitle

\section{Introduction}
\label{Introduction}
A main open problem in modern Galois theory is the characterization of the
profinite groups which are realizable as absolute Galois groups of fields $F$.
The torsion in such groups is described by the Artin--Schreier theory from the late 1920's,
namely, it consists solely of involutions.
More refined information on the structure of absolute Galois groups is given by Galois cohomology,
systematically developed starting the 1950's by Tate, Serre, and others.
Yet, explicit examples of torsion-free profinite groups which are not absolute Galois groups are rare.
In 1970, Milnor \cite{Milnor70} introduced his $K$-ring functor $K^M_*(F)$, and pointed out
close connections between this graded ring and the mod-$2$ Galois cohomology of the field.
This connection, in a more general form, became  known as the Bloch--Kato conjecture:
it says that for all $r\geq0$ and all $m$ prime to  $\Char\, F$,
there is a canonical isomorphism $K^M_r(F)/m\to H^r(G_F,\mu_m^{\tensor r})$ (\cite{GilleSzamuely}; see notation below).
The conjecture was proved for $r=2$ by Merkurjev and Suslin \cite{MerkurjevSuslin82},
for $r$ arbitrary and $m=2$ by Voevodsky \cite{Voevodsky03a},
and in general by Rost, Voevodsky, with a patch by Weibel (\cite{Voevodsky03b}, \cite{Weibel09},
\cite{Weibel08}, \cite{HaesemeyerWeibel09}).

In this paper we obtain new constrains on the group structure of
absolute Galois groups of fields, using this isomorphism.
We use these constrains to produce new examples of torsion-free profinite groups which are not absolute Galois groups.
We also demonstrate that the maximal pro-$p$ quotient of the absolute
Galois group can be characterized in purely cohomological terms (see Theorem \ref{thm 8.5}).
The main object of our paper is a remarkable small quotient of the absolute Galois group,
which, because of the above isomorphism, already carries a substantial information about the arithmetic of $F$.

More specifically, fix a prime number $p$ and a $p$-power $q=p^d$, with $d\geq1$.
All fields which appear in this paper will be tacitly assumed to contain a primitive $q$th root of unity.
Let $F$ be such a field and let $G_F=\Gal(F_\sep/F)$ be its absolute Galois group, where $F_\sep$ is
the separable closure of $F$.
Let $H^*(G_F)=H^*(G_F,\dbZ/q)$ be the Galois cohomology ring with the trivial action of $G_F$ on $\dbZ/q$.
Our new constraints relate the descending $q$-central sequence $G_F^{(i)}$,
$i=1,2,3\nek$ of $G_F$ (see \S\ref{section the descending central sequence})
with $H^*(G_F)$.
Setting $G_F^{[i]}=G_F/G_F^{(i)}$, we show that the quotient $G_F^{[3]}$ determines $H^*(G_F)$, and vice versa.
Specifically, we prove:

\begin{thmA}
The inflation map gives an isomorphism
\[
H^*(G_F^{[3]})_\dec\xrightarrow{\sim} H^*(G_F),
\]
where $H^*(G_F^{[3]})_\dec$ is the decomposable part of $H^*(G_F^{[3]})$
(i.e., its subring generated by degree $1$ elements).
\end{thmA}

We further have the following converse results.

\begin{thmB}
$G_F^{[3]}$ is uniquely determined by  $H^r(G_F)$  for $r=1,2$, the cup product
$\cup\colon H^1(G_F)\times H^1(G_F)\to H^2(G_F)$ and the Bockstein
homomorphism $\beta\colon H^1(G_F)\to H^2(G_F)$ (see \S2 for the definition of
$\beta$).
\end{thmB}

\begin{thmC}
Let $F_1$, $F_2$ be fields and let $\pi\colon G_{F_1}\to G_{F_2}$ be a
(continuous) homomorphism.
The following conditions are equivalent:
\begin{enumerate}
\item[(i)]
the induced map $\pi^*\colon H^*(G_{F_2})\to H^*(G_{F_1})$ is an isomorphism;
\item[(ii)]
the induced map $\pi^{[3]}\colon G_{F_1}^{[3]}\to G_{F_2}^{[3]}$ is an
isomorphism.
\end{enumerate}
\end{thmC}

Theorems A--C show that $G_F^{[3]}$ is a Galois-theoretic analog of the
cohomology ring $H^*(G_F)$.
Its structure is considerably simpler and more accessible than the full absolute Galois group $G_F$
(see e.g., \cite{EfratMinac}).
Yet, as shown in our theorems, these small and accessible quotients encode and control the entire
cohomology ring.
Results similar to Theorems A--C are valid in a relative pro-$p$ setting,
where one replaces $G_F$ by its maximal pro-$p$ quotient $G_F(p)=\Gal(F(p)/F)$
(here $F(p)$ is the compositum of all finite Galois extensions of $F$ of
$p$-power order; see Remark \ref{www}).
We further show:

\begin{thmD}
Let $F_1$, $F_2$ be fields containing a root of unity of order $p$ and let $\pi\colon G_{F_1}(p)\to G_{F_2}(p)$ be a
(continuous) homomorphism.
Then $\pi$ is an isomorphism if and only if the induced map $\pi^{[3]}\colon G_{F_1}^{[3]}\to G_{F_2}^{[3]}$ is an
isomorphism.
\end{thmD}

In the case $q=2$ the group $G_F^{[3]}$ has been extensively studied under the name ``$W$-group",
in particular in connection with quadratic forms (\cite{Spira87}, \cite{MinacSpira90}, \cite{MinacSpira96},
\cite{AdemKaraMinac99}, \cite{MaheMinacSmith04}). In this special case,
Theorem A was proved in \cite{AdemKaraMinac99}*{Th.\ 3.14}. It was further
shown that then $G_F^{[3]}$ has great arithmetical significance: it encodes
large parts of the arithmetical structure of $F$, such as its orderings, its
Witt ring, and certain non-trivial valuations. Theorem A explains this
surprising phenomena, as these arithmetical objects are known to be encoded in
$H^*(G_F)$ (with the additional knowledge of the Kummer element of $-1$).

First links between these quotients and the Bloch--Kato conjecture, and its special case the Merkurjev--Suslin theorem,
were already noticed in a joint work of Min\'a\v c and Spira in \cite{Spira87} and in Bogomolov's paper \cite{Bogomolov92}.
The latter paper was the first in a remarkable line of works by Bogomolov and Tschinkel
(\cite{Bogomolov92}, \cite{BogomolovTschinkel08}, \cite{BogomolovTschinkel09}, \cite{BogomolovTschinkel10}), as well as by Pop (unpublished),
focusing on the closely related quotient $G_F/[G_F,[G_F,G_F]]$
(the analog of $G_F^{[3]}$ for $q=0$), where $F$ is a function field over an algebraically closed field.
There the viewpoint is that of ``birational anabelian geometry":
namely,  it is shown that for certain important classes of such function fields, $F$ itself is determined by this quotient.
Our work, on the other hand, is aimed at clarifying the structure of the smaller Galois group $G_F^{[3]}$ and its
connections with the Galois cohomology and arithmetic of almost arbitrary fields,
focusing on the structure of absolute Galois groups.

Our approach is purely group-theoretic, and the main results above are in fact
proved for arbitrary profinite groups which satisfy certain conditions on
their cohomology (Theorem \ref{zzz}, Proposition \ref{Phi is G3}, Theorem \ref{equivalence}, and Remark \ref{pro-p}).
A key point is a rather general group-theoretic approach, partly inspired by \cite{GaoMinac97},
to the Milnor $K$-ring construction by means of quadratic hulls
of graded algebras (\S\ref{section graded rings}).
The Rost--Voevodsky theorem on the bijectivity of the Galois symbol shows that these
cohomological conditions are satisfied by absolute Galois groups as above.
Using this we deduce in \S\ref{section on examples} Theorems A--D in their field-theoretic version.

We wish to thank the referee for his insightful and thoughtful remarks and suggestions.
We also thank L.\ Avramov, M.\ Behrens, T.\ Chinburg, J.-L.\ Colliot-Th\'el\`ene, W.-D.\ Geyer, P.\ Goerss,
M.\ Hovey, M.\ Jarden, P.\ May, G.\ Prasad, Z.\  Reichstein and T.\ Szamuely for their comments related to talks on this work given
at the 2009 Field Arithmetic meeting in Oberwolfach, the AMS 2009 meeting in University Park, Pennsylvania, the University
of Chicago, Northwestern University, and the University of Nebraska, Lincoln.
We are also grateful to Thong Nguyen Quang Do for his interest in this work.
In particular, after our paper was posted on the arXiv, we were informed by him that
he also obtained a variant of Theorem \ref{RS}, and proved Theorem A in degree 2 using embedding problem techniques.
We thank him for sending us his (unpublished) notes \cite{NguyenQuangDo96}.

\section{Cohomological preliminaries
\label{section cohomological preliminaries}}
We work in the category of profinite groups.
Thus subgroups are always tacitly assumed to be closed and homomorphism are assumed to be continuous.
For basic facts on Galois cohomology we refer e.g., to \cite{NeukirchSchmidtWingberg}, \cite{SerreCG}, or \cite{Koch02}.
We abbreviate $H^r(G)=H^r(G,\dbZ/q)$ with the trivial $G$-action on $\dbZ/q$.
Let $H^*(G)=\bigoplus_{r=0}^\infty H^r(G)$ be the graded cohomology ring with the cup product.
We write $\res$, $\inf$, and $\trg$ for the restriction, inflation, and transgression maps, respectively.
Given a homomorphism $\pi\colon G_1\to G_2$ of profinite groups, we write $\pi^*\colon H^*(G_2)\to H^*(G_1)$ and
$\pi^*_r\colon H^r(G_2)\to H^r(G_1)$ for the induced homomorphisms.
The \textbf{Bockstein homomorphism} $\beta_G\colon H^1(G)\to H^2(G)$ of $G$ is the
connecting homomorphism arising from the short exact sequence of trivial $G$-modules
\[
0\ \to\ \dbZ/q\ \to\ \dbZ/q^2\ \to\ \dbZ/q\ \to\ 0.
\]
When $q=2$ one has $\beta_G(\psi)=\psi\cup\psi$ \cite{EfratMinac}*{Lemma 2.4}.

Given a normal subgroup $N$ of $G$, there is a natural action of $G$ on $H^r(N)$.
For $r=1$ it is given by $\varphi\mapsto\varphi^g$, where
$\varphi^g(n)=\varphi(g\inv ng)$ for $\varphi\in H^1(N)$, $g\in G$ and $n\in N$.
Let $H^r(N)^G$ be the group of all $G$-invariant elements of $H^r(N)$.
Recall that there is a $5$-term exact sequence
\[
0\to H^1(G/N)\xrightarrow{\inf_G} H^1(G)\xrightarrow{\res_N} H^1(N)^G\xrightarrow{\trg_{G/N}}
 H^2(G/N) \xrightarrow{\inf_G}H^2(G),
\]
which is functorial in $(G,N)$ in the natural sense.

\section{Graded rings\label{section graded rings}}
Let $R$ be a commutative ring and $\calA=\bigoplus_{r=0}^\infty A_r$ a graded
associative $R$-algebra with $A_0=R$.
Assume that $\calA$ is either commutative or graded-commutative
(i.e., $ab=(-1)^{rs}ba$ for $a\in A_r$, $b\in A_s$).
For $r\geq0$ let $A_{\dec,r}$ be the $R$-submodule of $A_r$ generated by all products of
$r$ elements of $A_1$ (by convention $A_{\dec,0}=R$).
The graded $R$-subalgebra $\calA_\dec=\bigoplus_{r=0}^\infty A_{\dec,r}$ is the
\textbf{decomposable part} of $\calA$.
We say that $A_r$ (resp., $\calA$) is \textbf{decomposable} if $A_r=A_{\dec,r}$ (resp., $\calA=\calA_\dec$).

Motivated by the Milnor $K$-theory of a field \cite{Milnor70},  we define the quadratic hull $\hat\calA$
of the algebra $\calA$ as follows.
For $r\geq0$ let $T_r$ be the $R$-submodule of $A_1^{\tensor r}$ generated by
all tensors $a_1\tensor\cdots\tensor a_r$ such that $a_ia_j=0\in A_2$
for some distinct $1\leq i,j\leq r$ (by convention, $A_1^{\tensor0}=R$, $T_0=0$).
We define $\hat\calA$ to be the graded $R$-algebra
$\hat\calA=\bigoplus_{r=0}^\infty A_1^{\tensor r}/T_r$ with
multiplicative structure induced by the tensor product.
Because of the commutativity/graded-commutativity, there is a canonical graded $R$-algebra epimorphism
$\omega_\calA\colon\hat\calA\to\calA_\dec$, which is the identity map in degree $1$.
We call $\calA$ \textbf{quadratic} if $\omega_\calA$ is an isomorphism.
Note that
\[
\hat\calA=(\hat\calA)_\dec=(\widehat{\calA_\dec}).
\]

These constructions are functorial in the sense that every graded $R$-algebra
morphism $\varphi=(\varphi_r)_{r=0}^\infty\colon \calA\to\calB$ induces in a
natural way graded $R$-algebra morphisms
\[
\varphi_\dec=(\varphi_{\dec,r})_{r=0}^\infty\colon\calA_\dec\to\calB_\dec,\qquad
\hat\varphi=(\hat\varphi_r)_{r=0}^\infty\colon \hat \calA\to\hat\calB
\]
with a commutative square
\begin{equation}
\label{cd} \xymatrix{
\hat\calA \ar[r]^{\hat\varphi} \ar[d]_{\omega_\calA} &  \hat\calB \ar[d]^{\omega_\calB}\\
\calA_\dec \ar[r]^{\varphi_\dec} & \ \calB_\dec. }
\end{equation}

The proof of the next fact is straightforward.

\begin{lem}
\label{bijectivity of hat varphi}
$\hat\varphi$ is an isomorphism if and only
if $\varphi_1$ is an isomorphism and $\varphi_{\dec,2}$ is a monomorphism.
\end{lem}

\begin{rem}
\rm
When $G$ is a profinite group, $R=\dbZ/2$, and  $\calA=H^*(G,\dbZ/2)$ the ring $\hat \calA$ coincides with
the ring $\mathrm{Mil}(G)$ introduced and studied in \cite{GaoMinac97}.
In the case where $G=G_F$ for a field $F$ as before, this ring is naturally isomorphic to $K^M_*(F)/2$.
Thus in this way one can construct $K^M_*(F)/p$ for any $p$ in a purely group-theoretic way.
\end{rem}

\section{The descending central sequence\label{section the descending central sequence}}
Let $G$ be a profinite group and let $q=p^d$ be either a $p$-power or $0$.
The \textbf{descending $q$-central sequence} of $G$ is defined inductively by
\[
G^{(1,q)}=G,\quad G^{(i+1,q)}=(G^{(i,q)})^q[G^{(i,q)},G],\quad  i=1,2,\ldots\
.
\]
Thus $G^{(i+1,q)}$ is the closed subgroup of $G$ topologically generated by
all powers $h^q$ and all commutators $[h,g]=h\inv g\inv hg$, where $h\in G^{(i,q)}$ and $g\in G$.
Note that $G^{(i,q)}$ is normal in $G$.
For $i\geq1$ let $G^{[i,q]}=G/G^{(i,q)}$.
When $q=0$ the sequence $G^{(i,0)}$ is called the \textbf{descending central sequence} of $G$.
Usually $q$ will be fixed, and  we will abbreviate
\[
G^{(i)}=G^{(i,q)},\ G^{[i]}=G^{[i,q]}.
\]
We will be mostly interested in these groups for $i=1,2,3$.
Note that $G^{[1]}=1$, and $G^{[2]}$ is the maximal quotient of $G$ of the form $\prod_{l=1}^d(\dbZ/p^l)^{I_l}$
for some index sets $I_1\nek I_d$.
In the case $i=3$ we have

\begin{lem}
\label{extension}
For $G_\ab=G/[G,G]$ one has a central extension
\[
1\to \bigl([G,G]/[[G,G],G]([G,G]\cap G^{q^2})\bigr)/q\to G^{[3]}\to G_\ab/q^2\to1.
\]
\end{lem}
\begin{proof}
There is an exact sequence
\[
1\to [G,G]/([G,G]\cap G^{(3)}) \to G^{[3]} \to G_\ab/q^2\to 1 .
\]
One has the identities (see \cite{Labute66}*{Prop.\ 5} and its proof)
\[
(g_1g_2)^q\equiv g_1^qg_2^q[g_2,g_1]^{\binom q2}, \quad
[g_1g_2,h]\equiv[g_1,h][g_2,h] \pmod{[[G,G],G]}
\]
It follows that $[G^{(2)},G]\leq [G,G]^q[[G,G],G]$, from which we get
\[
[G,G]\cap G^{(3)}=[[G,G],G]([G,G]\cap G^{q^2})[G,G]^q.
\qedhere
\]
\end{proof}

Any profinite  homomorphism (resp., epimorphism) $\pi\colon G\to H$ restricts
to a homomorphism (resp.,  an epimorphism)  $\pi^{(i)}\colon G^{(i)}\to H^{(i)}$.
Hence $\pi$ induces a homomorphism (resp., an epimorphism) $\pi^{[i]}\colon G^{[i]}\to H^{[i]}$.

\begin{lem}
\label{yyy}
Let $\pi\colon G_1\to G_2$ be an epimorphism of profinite groups.
Then $\Ker(\pi^{[i]})=\Ker(\pi)G_1^{(i)}/G_1^{(i)}$ for all $i\geq1$.
\end{lem}
\begin{proof}
The map $G_1\to G_2^{[i]}=\pi(G_1)/\pi(G_1)^{(i)}$ induced by $\pi$ has kernel
$\Ker(\pi)G_1^{(i)}$, whence the assertion.
\end{proof}

We will also need the following result of Labute \cite{Labute66}*{Prop.\ 1 and 2}
(see also \cite{NeukirchSchmidtWingberg}*{Prop.\ 3.9.13}).

\begin{prop}
\label{uniqueness of presentation}
Let $S$ be a free pro-$p$ group on generators $\sig_1\nek\sig_n$.
Consider the Lie $\dbZ_p$-algebra $\gr(S)=\bigoplus_{i=1}^\infty S^{(i,0)}/S^{(i+1,0)}$,
with Lie brackets induced by the commutator map.
Then $\gr(S)$ is a free Lie $\dbZ_p$-algebra on the images of $\sig_1\nek\sig_n$ in $\gr_1(S)$.
In particular, $S^{(2,0)}/S^{(3,0)}$ has a system of representatives consisting of all products
$\prod_{1\leq k<l\leq n}[\sig_k,\sig_l]^{a_{kl}}$, where $a_{kl}\in\dbZ_p$.
\end{prop}

\section{Duality\label{section duality}}

From now on assume that $q=p^d$ with $p$ prime. 
One has

\begin{prop}[\cite{EfratMinac}*{Cor.\ 2.2}]
\label{duality}
For a normal subgroup $R$ of a profinite group $G$, there is a perfect duality
\[
R/R^q[R,G]\ \times\ H^1(R)^G\ \to\ \dbZ/q ,\qquad (\bar r,\varphi)\mapsto
\varphi(r).
\]
\end{prop}

We observe that this duality is functorial in the following sense:

\begin{prop}
\label{functoriality} For a homomorphism $\pi\colon G_1\to G_2$ of profinite
groups and for normal subgroups $R_i$ of $G_i$, $i=1,2$, such that
$\pi(R_1)\leq R_2$, the following induced diagram commutes in the natural sense:
\[
\xymatrix{
R_1/R_1^q[R_1,G_1] \ar[d]_{\overline{\pi|_{R_1}}} & *-<3pc>{\times} & H^1(R_1)^{G_1}\ar[r] &  \dbZ/q\ar@{=}[d] \\
R_2/R_2^q[R_2,G_2]  & *-<3pc>{\times} & H^1(R_2)^{G_2}\ar[r]\ar[u]_{(\pi|_{R_1})_1^*} &  \dbZ/q. \\
}
\]
\end{prop}

We list some consequences of this duality.

\begin{lem}
\label{induced map 2} Let $\pi\colon G_1\to G_2$ be a homomorphism of profinite groups.
Then $\pi^{[2]}$ is an isomorphism if and only if $\pi^*_1\colon H^1(G_2)\to H^1(G_1)$ is an isomorphism.
In this case, $\pi^*_\dec\colon H^*(G_2)_\dec\to H^*(G_1)_\dec$ is surjective.
\end{lem}
\begin{proof}
For the equivalence use Proposition \ref{functoriality} with $G_i=R_i$,
$i=1,2$. The second assertion follows from the first one.
\end{proof}

\begin{lem}
\label{triviality of restriction}
$\res\colon H^1(G^{(i)})^G\to H^1(G^{(i+1)})^G$ is trivial for all $i\geq1$.
\end{lem}
\begin{proof}
The homomorphism $G^{(i+1)}/G^{(i+2)}\to G^{(i)}/G^{(i+1)}$ induced by the
inclusion map $G^{(i+1)}\hookrightarrow G^{(i)}$ is trivial.
Now use Proposition \ref{functoriality} with $G_1=G_2=G$ , $R_1=G^{(i+1)}$, and $R_2=G^{(i)}$.
\end{proof}

The next proposition is a variant of \cite{Bogomolov92}*{Lemma 3.3}.

\begin{lem}
\label{equal kernels} The inflation maps
\[
\textstyle\inf_{G^{[3]}}\colon H^2(G^{[2]})\to H^2(G^{[3]}),\quad
\textstyle\inf_G\colon H^2(G^{[2]})\to H^2(G)
\]
have the same kernel.
\end{lem}
\begin{proof}
The functoriality of the $5$-term sequence (see \S\ref{section cohomological
preliminaries} and \cite{EfratMinac}*{\S2B)}) gives a commutative diagram with
exact rows
\[
\xymatrix{
  H^1(G^{(2)})^G \ar[r]^{\trg_{G^{[2]}}}\ar[d]_{\res_{G^{(3)}}}  & H^2(G^{[2]}) \ar[r]^{\ \inf_G}\ar[d]_{\inf_{G^{[3]}}}
     &  H^2(G)\ar@{=}[d] \\
  H^1(G^{(3)})^G \ar[r]^{\trg_{G^{[3]}}}  &  H^2(G^{[3]}) \ar[r]^{\ \inf_G} & H^2(G).
}
\]
By Lemma \ref{triviality of restriction}, $\res_{G^{(3)}}=0$, and the assertion follows by a diagram chase.
\end{proof}

\begin{prop}
\label{image of inflation is decomposable}
Let $G$ be a profinite group with $H^2(G)$ decomposable.
Then $\inf_{G^{[3]}}$ maps $H^2(G^{[2]})$ into $H^2(G^{[3]})_\dec$.
\end{prop}
\begin{proof}
Since $\inf_G\colon H^1(G^{[2]})\to H^1(G)$ is an isomorphism,
the map $\inf_G\colon H^2(G^{[2]})_\dec\to H^2(G)_\dec=H^2(G)$ is surjective.
Furthermore, $\inf_{G^{[3]}}$ maps $H^2(G^{[2]})_\dec$ into $H^2(G^{[3]})_\dec$.
Now apply Lemma 5.5.
\end{proof}

\begin{prop}
\label{injectivity of functor}
Let $\pi\colon G_1\to G_2$ be a homomorphism of profinite groups with $\pi_1^*$ bijective,
and with $\pi^*_2$ injective on the image $A$ of $\inf_{G_2}\colon H^2(G_2^{[2]})\to H^2(G_2)$.
Then $\pi^{[3]}$ is an isomorphism.
\end{prop}
\begin{proof}
By Lemma \ref{induced map 2}, the maps $\pi^{[2]}\colon G_1^{[2]}\to
G_2^{[2]}$ and $\inf_{G_i}\colon H^1(G_i^{[2]})\to H^1(G_i)$, $i=1,2$, are isomorphisms.
By the functoriality of the $5$-term sequence, there is a commutative diagram with exact rows
\[
\xymatrix{ 0\ar[r] & H^1(G_2^{(2)})^{G_2} \ar[d]^{(\pi^{(2)})_1^*}
\ar[r]^{\quad\trg} &
H^2(G_2^{[2]}) \ar[r]^{\inf_{G_2}}\ar[d]^{(\pi^{[2]})^*_2} & A\ar[r]\ar[d]^{\pi_2^*} & 0 \\
0\ar[r] & H^1(G_1^{(2)})^{G_1} \ar[r]^{\quad\trg} & H^2(G_1^{[2]})
\ar[r]^{\inf_{G_1}} & H^2(G_1). }
\]

Now since $\pi^{[2]}$ is an isomorphism, so is $(\pi^{[2]})^*_2$.
By assumption, $\pi_2^*$ is injective on $A$.
A snake lemma argument shows that $(\pi^{(2)})_1^*$ is an isomorphism.
By passing to duals (using Proposition \ref{functoriality} with $R_i=G_i^{(2)}$), we obtain that the map
$\bar\pi\colon G_1^{(2)}/G_1^{(3)}\to G_2^{(2)}/G_2^{(3)}$ induced by $\pi$ is also an isomorphism.
Thus in the commutative diagram
\[
\xymatrix{
1\ar[r] & G_1^{(2)}/G_1^{(3)} \ar[r]\ar[d]_{\bar\pi}  & G_1^{[3]} \ar[r]\ar[d]^{\pi^{[3]}}
& G_1^{[2]}\ar[d]^{\pi^{[2]}}\ar[r]  & 1 \\
1\ar[r] & G_2^{(2)}/G_2^{(3)}\ar[r]  & G_2^{[3]}\ar[r]  & G_2^{[2]}  \ar[r] &
1 }
\]
both $\bar\pi$ and $\pi^{[2]}$ are isomorphisms.
By the snake lemma, so is $\pi^{[3]}$.
\end{proof}

\section{Morphisms of cohomology rings\label{section Morphisms of cohomology rings}}
We will now use the general constructions of \S\ref{section graded rings} with
the base ring $R=\dbZ/q$.

\begin{lem}
\label{pidec isomorphism}
Let $\pi\colon G_1\to G_2$ be a homomorphism of
profinite groups such that $\widehat{\pi^*}$ is an isomorphism and  $H^*(G_1)$ is quadratic.
Then $\pi_\dec^*$ is an isomorphism and $H^*(G_2)$ is also quadratic.
\end{lem}
\begin{proof}
By (\ref{cd}), there is a commutative square
\[
\xymatrix{
\widehat{H^*(G_2)} \ar[d]_{\omega_{H^*(G_2)}} \ar[r]^{\widehat{\pi^*}} & \widehat{H^*(G_1)}\ar[d]^{\omega_{H^*(G_1)}} \\
H^*(G_2)_\dec \ar[r]^{\pi^*_\dec}  &  H^*(G_1)_\dec. }
\]
Now use the surjectivity of $\omega_{H^*(G_2)}$ and the assumptions.
\end{proof}

\begin{prop}
\label{induced map} Let $\pi\colon G_1\to G_2$ be a homomorphism of profinite groups such that
$\pi^{[3]}$ is an isomorphism.
Then:
\begin{enumerate}
\item[(a)]
$\pi^*_{\dec,2}\colon H^2(G_2)_\dec\to H^2(G_1)_\dec$ is an isomorphism;
\item[(b)]
$\widehat{\pi^*}\colon\widehat{H^*(G_2)}\to \widehat{H^*(G_1)}$ is an
isomorphism.
\end{enumerate}
\end{prop}
\begin{proof}
First note that since $\pi^{[3]}$ is an isomorphism, so is $\pi^{[2]}$.

\medskip

(a)\quad One has a commutative diagram with isomorphisms as indicated:
\[
\xymatrix{ H^2(G_2^{[2]})_\dec \ar[r]^{\inf}
\ar[d]_{\wr}^{(\pi^{[2]})_{\dec,2}^*} & H^2(G_2^{[3]})_\dec
\ar[d]_{\wr}^{(\pi^{[3]})_{\dec,2}^*} \ar[r]^{\inf} &
H^2(G_2)_\dec\ar[d]^{\pi_{\dec,2}^*}\\
H^2(G_1^{[2]})_\dec \ar[r]^{\inf}  & H^2(G_1^{[3]})_\dec  \ar[r]^{\inf} &
H^2(G_1)_\dec. }
\]
Further, by Lemma \ref{induced map 2}, $\inf_{G_i}\colon
H^2(G_i^{[2]})_\dec\to H^2(G_i)_\dec$ is surjective for $i=1,2$.
The assertion now follows using Lemma \ref{equal kernels} (for $G=G_1$).

\medskip

(b)\quad
By Lemma \ref{induced map 2}, $\pi_1^*\colon H^1(G_2)\to H^1(G_1)$ is
an isomorphism. Now use (a) and Lemma \ref{bijectivity of hat varphi}.
\end{proof}

Combining the previous results we obtain:

\begin{thm}
\label{equivalence}
Let $\pi\colon G_1\to G_2$ be a homomorphism of profinite
groups with $H^*(G_1)$ quadratic and $H^2(G_2)$  decomposable.
The following conditions are equivalent:
\begin{enumerate}
\item[(a)]
$\pi_1^*$ is bijective and $\pi^*_2$ is injective;
\item[(b)]
$\pi^{[3]}$ is an isomorphism;
\item[(c)]
$\widehat{\pi^*}$ is an isomorphism;
\item[(d)]
$\pi^*_\dec$ is an isomorphism.
\end{enumerate}
\end{thm}
\begin{proof}
(a)$\Rightarrow$(b) follows from Proposition \ref{injectivity of functor}.
(b)$\Rightarrow$(c) is Proposition \ref{induced map}(b).
For (c)$\Rightarrow$(d) use Lemma \ref{pidec isomorphism}.
Finally, (d)$\Rightarrow$(a) follows from the decomposability of $H^2(G_2)$.
\end{proof}

\begin{rem}
\label{pro-p}
\rm
When $G_1,G_2$ are pro-$p$ groups, condition (a) of Theorem \ref{equivalence} (and therefore all other conditions)
is equivalent to $\pi$ being an isomorphism \cite{Serre65}*{Lemma 2}.
\end{rem}

Denote the maximal pro-$p$ quotient of $G$ by $G(p)$.

\begin{lem}
\label{qqq} $\inf_G\colon H^2(G(p))\to H^2(G)$ is injective.
\end{lem}
\begin{proof}
The kernel $N$ of the epimorphism $G\to G(p)$ satisfies $H^1(N)^G=0$.
Now use the $5$-term sequence.
\end{proof}

\begin{thm}
\label{zzz}
Let $G$ be a profinite group with $H^*(G)$ quadratic.
Let $N$ be a normal subgroup of $G$ with $N\leq G^{(3)}$.
Then
\begin{enumerate}
\item[(a)]
$\inf_G\colon H^*(G/N)_\dec\to H^*(G)_\dec$ is an isomorphism.
\item[(b)]
$H^*(G/N)$ is quadratic.
\item[(c)]
If $H^2(G/N)$ is decomposable, then $N$ is contained in the kernel of the
canonical map $G\to G(p)$.
\end{enumerate}
\end{thm}
\begin{proof}
Consider the natural epimorphism $\pi\colon G\to G/N$.
By Lemma \ref{yyy}, $\pi^{[3]}$ is an isomorphism.
By Proposition \ref{induced map}(b), $\widehat{\pi^*}$ is also an isomorphism.
Assertions (a) and (b) now follow from Lemma \ref{pidec isomorphism}.

To prove (c), let $M$ be the kernel of the projection $G\to G(p)$ and set $R=MN$.
A standard group-theoretic argument shows that $H^1(M)=0$.
Since $H^1(R/N)\isom H^1(M/M\cap N)$ injects into it by inflation, it is also trivial.
By the 5-term sequence, $\res_N\colon H^1(R)\to H^1(N)$ is therefore injective.

By (a),  $\inf_G\colon H^2(G/N)=H^2(G/N)_\dec\to H^2(G)_\dec$ is an isomorphism.
Also, since $N,M\leq G^{(2)}$, the inflation maps $H^1(G/R)\to
H^1(G)$, $H^1(G/N)\to H^1(G)$ are isomorphisms (Lemma \ref{induced map 2}).
Using as before the functoriality of the $5$-term sequence, we get a
commutative diagram with exact rows
\[
\xymatrix{ 0 \ar[r] & H^1(R)^G\ar[r]^{\trg_{G/R}} \ar@{>->}[d]_{\res_N} &
H^2(G/R)
\ar[r]^{\inf_G}\ar[d]^{\inf_{G/N}} &  H^2(G) \\
0 \ar[r] & H^1(N)^G\ar[r]^{\trg_{G/N}} & H^2(G/N) \ar[r]^{\inf_G}_{\sim} &
H^2(G)_\dec \ar@{^{(}->}[u]. }
\]
Consequently, $H^1(R)^G=0$, whence  $H^1(R/M,\dbZ/p)^{G(p)}=0$.
Since $G(p)$ is pro-$p$, it follows as in \cite{Serre65}*{Lemma 2} that $R/M=1$, i.e., $N\leq M$.
\end{proof}

\section{Cohomology determines $G^{[3]}$ \label{cohomology detrmines G3}}

Let $G$ be a profinite group, $S$ a free pro-$p$ group, and $\pi\colon S\to G(p)$ an epimorphism.
Let $\varphi\colon G(p)\to G^{[2]}$ be the natural map and set $R=\Ker(\pi)$ and $T=\Ker(\varphi\circ\pi)$.
One has a commutative diagram with exact rows
\[
\xymatrix{
1\ar[r] & R\ar[r]\ar@{_{(}->}[d]  & S\ar[r]^{\pi}\ar@{=}[d] & G(p)\ar[r]\ar[d]^{\varphi} & 1 \\
1\ar[r] & T\ar[r] & S\ar[r] & G^{[2]} \ar[r] & 1.
}
\]
As $H^2(S)=0$, the corresponding $5$-term sequences give as before a
commutative diagram with exact rows
\begin{equation}
\label{cd3}
\xymatrix{
0\ar[r] & H^1(S)/H^1(G(p)) \ar[r] & H^1(R)^S\ar[r]^{\trg} & H^2(G(p))\ar[r] & 0 \\
0\ar[r] &  H^1(S)/H^1(G^{[2]})\ar[r]\ar@{->>}[u] & H^1(T)^S\ar[r]^{\trg}\ar[u]^{\res_R} &
 H^2(G^{[2]})\ar[u]_{\inf_G}\ar[r] & 0.}
\end{equation}
Proposition \ref{duality} gives a commutative diagram of perfect pairings
\begin{equation}
\label{cd4}
\xymatrix{
R/R^q[R,S] \ar[d]_{\iota} & *-<3pc>{\times} & H^1(R)^S \ar[r]  & \dbZ/q \ar@{=}[d]\\
T/T^q[T,S] & *-<3pc>{\times} & H^1(T)^S \ar[u]_{\res_R} \ar[r] &
\dbZ/q
. \\
}
\end{equation}

\begin{thm}
\label{group theoretic MS}
The following conditions are equivalent:
\begin{enumerate}
\item[(a)]
$\inf_G\colon H^2(G^{[2]})\to H^2(G(p))$ is surjective;
\item[(b)]
$\res_R\colon H^1(T)^S\to H^1(R)^S$ is surjective;
\item[(c)]
$R/R^q[R,S]\to T/T^q[T,S]$ is injective;
\item[(d)]
$R^q[R,S]=R\cap T^q[T,S]$.
\end{enumerate}
\end{thm}
\begin{proof}
(a)$\Leftrightarrow$(b): \quad
Apply the snake lemma to (\ref{cd3}) to obtain that
$\Coker(\res_R)\isom\Coker(\inf_G)$, and the equivalence follows.

\medskip

(b)$\Leftrightarrow$(c): \quad
Use (\ref{cd4}).

\medskip

(c)$\Leftrightarrow$(d): \quad
The homomorphism in (c) breaks as
\[
R/R^q[R,S] \to R/(R\cap T^q[T,S]) \to T/T^q[T,S],
\]
where the right map is injective.
Therefore the injectivity of the composed map is equivalent to that of the left map, i.e., to (d).
\end{proof}

Of special importance is the case where the presentation $1\to R\to S\xrightarrow{\pi} G\to 1$ is \text{\bf minimal},
i.e., $R\leq S^{(2)}$.
Then $\pi^{[2]}$ is an isomorphism, by Lemma \ref{yyy}.
Since $\varphi^{[2]}$ is always an isomorphism, $T=S^{(2)}$, so $T^q[T,S]=S^{(3)}$.
In this case, and assuming further that $q=2$, the equivalences of Theorem \ref{group theoretic MS}
were obtained in \cite{GaoMinac97}*{Th.\ 2} and  \cite{MinacSpira96}*{\S5}.

\begin{rem}
\label{vvv}
\rm
By Lemma \ref{induced map 2}, $\inf_G\colon H^2(G^{[2]})_\dec\to H^2(G(p))_\dec$ is surjective.
Thus, if $H^2(G(p))$ is decomposable, then condition (a), whence all other conditions of
Theorem \ref{group theoretic MS}, are satisfied.
This was earlier observed by T.\ W\"urfel \cite{Wurfel85}*{Prop.\ 8} for condition (d), in the case of
minimal presentations and the descending ($0$-)\,central sequence, and
under the additional assumption that the maximal abelian quotient $G_\ab$ is torsion-free.
\end{rem}

Next let $f$ be the composed homomorphism
\[
H^1(T)^S\ \xrightarrow{\trg}\  H^2(G^{[2]}) \ \xrightarrow{\inf_{G(p)}}\ H^2(G(p))\ \xrightarrow{\inf_G}\ H^2(G)\ .
\]
Let $\Ker(f)^\vee$ be the dual of $\Ker(f)$ in $T/T^q[T,S]$ with respect to the pairing in (\ref{cd4}).

\begin{prop}
\label{Phi is G3}
Assume that the presentation $1\to R\to S\to G\to1$ is minimal.
Then $G^{[3]}\isom S^{[3]}/\Ker(f)^\vee$.
\end{prop}
\begin{proof}
Since $R\leq S^{(2)}=T$, Proposition \ref{duality} shows that the restrictions
$H^1(S)\to H^1(R)^S$ and $H^1(S)\to H^1(T)^S$ are the zero maps.
Hence the transgression maps in (\ref{cd3}) are isomorphisms.
By the injectivity of $\inf_G\colon H^1(G(p))\to H^1(G)$ (Lemma \ref{qqq}), $\Ker(f)=\Ker(\res_R)$.
Therefore, by (\ref{cd4}),
\[
\Ker(f)^\vee=\Ker(\res_R)^\vee=\Img(\iota)=RS^{(3)}/S^{(3)}.
\]
In view of Lemma \ref{yyy}, $RS^{(3)}/S^{(3)}$ is the kernel of the induced map
$S^{[3]}\to (S/R)^{[3]}\isom G(p)^{[3]}=G^{[3]}$.
Consequently, $S^{[3]}/\Ker(f)^\vee\isom G^{[3]}$.
\end{proof}

\begin{rem}
\label{detrmined by cohomology}
\rm
Assume that the profinite group $G$ satisfies
\begin{enumerate}
\item[($*$)]
$G^{[2]}\isom (\dbZ/q)^I$ for some set $I$.
\end{enumerate}
Then one can find $S$ and $\pi$ as above such that the presentation is minimal.
Namely, take a free pro-$p$ group $S$ such that  $S^{[2]}\isom G^{[2]}$.
This lifts to a homomorphism  $S\to G$, which is surjective, by a Frattini argument.

We further claim that under assumption ($*$), $f$ is determined
by $H^r(G)$, $r=1,2$, the cup product $\cup\colon H^1(G)\tensor_\dbZ H^1(G)\to
H^2(G)$ and the Bockstein map $\beta_G\colon H^1(G)\to H^2(G)$.
Indeed, $H^1(G)$ determines $G^{[2]}$ which is its Pontrjagin dual.
Set
\[
\Omega(G)=(H^1(G)\tensor_\dbZ H^1(G))\oplus H^1(G)
\]
and consider the homomorphism
\[
\Lambda_G\colon \Omega(G)\to H^2(G), \quad
(\alp_1,\alp_2)\mapsto\cup\alp_1+\beta_G(\alp_2).
\]
Then $\Lambda\colon \Omega\to H^2$ is a natural transformation of contravariant functors.
For $\bar G\isom(\dbZ/q)^I$, the map $\Lambda_{\bar G}$ is surjective, by \cite{EfratMinac}*{Cor.\ 2.11}
(when $\bar G$ is not finitely
generated apply a limit argument, as in \cite{EfratMinac}*{Lemma 4.1}).
Now the canonical map $G\to G^{[2]}$ induces a commutative square
\[
\xymatrix{ \Omega(G^{[2]}) \ar[r]^{\sim}\ar[d]_{\Lambda_{G^{[2]}}}   &
\Omega(G)\ar[d]^{\Lambda_G} \\
H^2(G^{[2]})\ar[r]^{\inf_G} & H^2(G). }
\]
From the surjectivity of $\Lambda_{G^{[2]}}$ we see that $\inf_G$ is determined
by $\Lam_G$ (modulo the identification $H^1(G^{[2]})\isom H^1(G)$).
Hence $\inf_G$ is determined by $H^1(G)$, $H^2(G)$, $\cup$ and $\beta_G$.
Also, $S$ is determined by $H^1(G)$ only.
\end{rem}

\section{Absolute Galois groups
\label{section on absolute Galois groups}}
Let $F$ be a field of characteristic $\neq p$ (containing as always the group $\mu_q$ of $q$th roots of unity).
By fixing a primitive $q$th root of unity we may identify $\mu_q=\dbZ/q$ as $G_F$-modules.
Let $K^M_*(F)$ be again the Milnor $K$-ring of $F$ \cite{Milnor70}.
Define graded $\dbZ/q$-algebras $\calA=K^M_*(F)/q$ and $\calB=H^*(G_F)$.
Observe that they are graded-commutative and $\calA=\calA_\dec=\hat\calA$.
It can be shown that the Kummer isomorphism $F^\times/(F^\times)^q\isom H^1(G_F)$, $a(F^\times)^q\mapsto
(a)_F$, extends to an epimorphism $h_F\colon \calA\to \calB_\dec$, called the \textbf{Galois symbol epimorphism}
\cite{GilleSzamuely}*{\S4.6}.
It is bijective in degrees $1$ and $2$, by the Kummer isomorphism and the injectivity part of the
Merkurjev--Suslin theorem, respectively \cite{GilleSzamuely}.
By Lemma \ref{bijectivity of hat varphi}, $\hat h_F\colon\calA=\hat\calA\to\hat \calB$ is an isomorphism.
Furthermore, by the Rost--Voevodsky theorem (\cite{Voevodsky03b}, \cite{Weibel09}, \cite{Weibel08}),
$h_F$ itself is an isomorphism and
$\calB=\calB_\dec$.
By (\ref{cd}), $\omega_\calB$ is also an isomorphism, i.e., $\calB=H^*(G_F)$ is quadratic.

Theorem \ref{zzz}(a)(b) give:

\begin{cor}
For a field $F$ as above one has:
\begin{enumerate}
\item[(a)]
$\inf_{G_F}\colon H^*(G_F^{[3]})_\dec\to H^*(G_F)$ is an isomorphism;
\item[(b)]
$H^*(G_F^{[3]})$ is quadratic.
\end{enumerate}
\end{cor}

This establishes Theorem A, and shows that for $F$ as above, $G_F^{[3]}$ determines $H^*(G_F)$.

Theorem C follows from Theorem \ref{equivalence}.

Finally, the projection $F^\times/(F^\times)^q\to F^\times/(F^\times)^p$ is
obviously surjective. By the Kummer isomorphism, so is the functorial map
$H^1(G_F)=H^1(G_F,\dbZ/q)\to H^1(G_F,\dbZ/p)$.
Consequently, ($*$) of Remark \ref{detrmined by cohomology} is satisfied for $G=G_F$.
Proposition \ref{Phi is G3} and Remark \ref{detrmined by cohomology} now give Theorem B.

\begin{rem}
\label{www}
\rm
For the maximal pro-$p$ Galois extension $F(p)$ of $F$ one has $H^1(G_{F(p)})=0$.
Hence $H^*(G_{F(p)})=H^*(G_{F(p)})_\dec$ vanishes in positive degrees.
A spectral sequence argument therefore shows that Lyndon--Hochschild--Serre spectral sequence
corresponding to $G_F\to G_F(p)$ collapses at $E_2$, giving an isomorphism $\inf\colon H^*(G_F(p))\to H^*(G_F)$
\cite{NeukirchSchmidtWingberg}*{Lemma 2.1.2}.
Therefore in the previous discussion and in Theorems A--C we may replace $G_F$ by $G_F(p)$.
In particular, $H^*(G_F(p))$ is quadratic and decomposable.
\end{rem}

In view of this remark, Theorem D follows from Remark \ref{pro-p} with $G_i=G_{F_i}(p)$, $i=1,2$.

We have the following variant of Theorem C for homomorphisms of the groups $G_{F_i}^{[3]}$.

\begin{thm}
\label{G3 variant of theorem C}
Let $F_1,F_2$ be fields as above
and let $\bar\pi\colon G_{F_1}^{[3]}\to G_{F_2}^{[3]}$ be a homomorphism.
The following conditions are equivalent:
\begin{enumerate}
\item[(a)]
$\bar\pi\colon G_{F_1}^{[3]}\to G_{F_2}^{[3]}$ is an isomorphism;
\item[(b)]
$\pi^*_\dec\colon H^*(G_{F_2}^{[3]})_\dec\to H^*(G_{F_1}^{[3]})_\dec$ is an isomorphism;
\item[(c)]
$\bar\pi^*_1\colon H^1(G_{F_2}^{[3]})\to H^1(G_{F_1}^{[3]})$ is bijective and
$\bar\pi^*_{2,\dec}\colon H^2(G_{F_2}^{[3]})_\dec\to H^2(G_{F_1}^{[3]})_\dec$ is injective.
\end{enumerate}
\end{thm}
\begin{proof}
(a)$\Rightarrow$(b)$\Rightarrow$(c): \quad
Trivial.

\medskip

(c)$\Rightarrow$(a): \quad
Since $H^2(G_{F_2})$ is decomposable, so is the image of
$\inf_{G_{F_2}^{[3]}}\colon H^2(G_{F_2}^{[2]})\to H^2(G_{F_2}^{[3]})$,
by Proposition \ref{image of inflation is decomposable}.
By assumption, $\bar\pi^*_2$ is therefore injective on this image.
Now apply Proposition \ref{injectivity of functor} with $G_i=G_{F_i}^{[3]}$, $i=1,2$.
\end{proof}

By Remarks \ref{vvv} and \ref{www}, the equivalent conditions of Theorem
\ref{group theoretic MS} are satisfied for $G_F$.
We therefore get:

\begin{thm}
\label{RS}
Let $1\to R\to S\to G_F(p)\to 1$ be a minimal presentation of $G_F(p)$ using
generators and relations, where $S$ is a free pro-$p$ group.
Then $R^q[R,S]=R\cap S^{(3)}$.
\end{thm}

Moreover, by Theorem \ref{zzz}(c), $G_F(p)$ has the following minimality
property with respect to decomposability of $H^2$.
Here $F^{(3)}$ is the fixed field of $G_F^{(3)}$ in $F_\sep$.

\begin{thm}
\label{thm 8.5}
Suppose that $K$ is a Galois extension of $F$ containing $F^{(3)}$.
Assume that $H^2(\Gal(K/F))$ is decomposable.
Then $K\supseteq F(p)$.
\end{thm}

This can be viewed as a cohomological characterization of the extension $F(p)/F$.

\section{Groups which are not maximal pro-$p$ Galois groups
\label{section on examples}}
We now apply Theorem A to give examples of pro-$p$ groups which cannot be realized
as maximal pro-$p$ Galois groups of fields (assumed as before to contain a root of unity of order $p$).
Our groups also cannot be realized as the \textsl{absolute} Galois group $G_F$ of any field $F$.
Indeed, assume this were the case.
When $\Char\,F=p$, the maximal pro-$p$ Galois group of $F$ is a free pro-$p$ group
\cite{NeukirchSchmidtWingberg}*{Th.\ 6.1.4},
whereas our groups are not free.
When $\Char\,F\neq p$, as $G_F$ is pro-$p$, necessarily $\mu_p\subseteq F$, and we get again a contradiction.

The groups we construct are only a sample of the most simple and straightforward examples illustrating our theorems, and
many other more complicated examples can be constructed along the same lines.

Throughout this section $q=p$.
We have the following immediate consequence of
the analog of Theorem A for maximal pro-$p$ Galois groups (Remark \ref{www}):

\begin{prop}
\label{principle}
If $G_1,G_2$ are pro-$p$ groups such that $G_1^{[3,p]}\isom G_2^{[3,p]}$ and
$H^*(G_1)\not\isom H^*(G_2)$, then at most one of them can be realized as the maximal pro-$p$ Galois group of a field.
\end{prop}

\begin{cor}
\label{cor to principle}
Let $S$ be a free pro-$p$ group and $R$ a nontrivial closed normal subgroup of $S^{(3,p)}$.
Then $G=S/R$ cannot occur as a maximal pro-$p$ Galois group of a field.
\end{cor}
\begin{proof}
Since $R\leq S^{(3,p)}\leq S^{(2,p)}$, there is an induced isomorphism $S^{[2,p]}\isom G^{[2,p]}$.
Hence $\inf\colon H^1(G)\to H^1(S)$ is an isomorphism.
As $H^2(S)=0$  \cite{NeukirchSchmidtWingberg}*{Cor. 3.9.5}, $\trg\colon H^1(R)^G\to H^2(G)$ is also an isomorphism.
Furthermore, $H^1(R)\neq0$, so a standard orbit counting argument shows that $H^1(R)^S\neq0$.
Hence $H^2(G)\neq0$.
Now $S\isom G_F$ for some field $F$ of characteristic $0$
\cite{FriedJarden}*{Cor.\ 23.1.2}, so we may apply Proposition \ref{principle} with $G_1=S$ and $G_2=G$.
\end{proof}

\begin{exam}  \rm
Let $S$ be a free pro-$p$ group on $2$ generators, and take $R=[S,[S,S]]$.
By Corollary \ref{cor to principle}, $G=S/R$ is not realizable as $G_F(p)$ for a field $F$ as above.
Note that $G/[G,G]\isom S/[S,S]\isom\dbZ_p^2$ and $[G,G]=[S,S]/S^{(3,0)}\isom\dbZ_p$ (Proposition \ref{uniqueness of
presentation}), so $G$ is torsion-free.
\end{exam}

\begin{exam}  \rm
Let $S$ be a free pro-$p$ group on $n\geq3$ generators $\sig_1\nek\sig_n$,  and let
$R$ be the closed normal subgroup of $S$ generated by $r=[\ldots[[\sig_1,\sig_2],\sig_3]\nek\sig_n]$.
By Corollary \ref{cor to principle}, $G=S/R$ is not realizable as $G_F(p)$ for a field $F$ as above.
By Proposition \ref{uniqueness of presentation}, $r\not\in (S^{(n,0)})^pS^{(n+1,0)}$.
Therefore \cite{Labute67b}*{Th.\ 4} implies that $G=S/R$ is torsion-free.
We further note that, by \cite{Romanovskij86}, $G$ contains a subgroup which is a free pro-$p$ group on two generators.
In particular, $G$ is not prosolvable.
\end{exam}

\begin{exam}  \rm
\label{example Labute-Gildenhuys} Let $S$ be a free pro-$p$ group on
generators $\sig_1,\sig_2$.
Let $r_1=\sig_1^p$ and $r_2=\sig_1^p[\sig_1,[\sig_1,\sig_2]]$.
Let $R_i$ be the normal subgroup of $S$ generated by $r_i$, and set $G_i=S/R_i$, $i=1,2$.
As $r_1r_2\inv\in S^{[3,p]}$ we have $G_1^{[3,p]}\isom G_2^{[3,p]}$.

Since $G_1$ has $p$-torsion, it has $p$-cohomological dimension
$\cd_p(G_1)=\infty$ \cite{SerreCG}*{I, \S3.3, Cor.\ 3}. On the other hand,
$\cd_p(G_2)=2$; indeed, this is proved in \cite{Labute67b}*{p.\ 144, Example
(1)} for $p\geq5$, in \cite{Labute67b}*{p.\ 157, Exemple} for $p=3$, and in
\cite{Guildenhuys68}*{Prop.\ 3} for $p=2$.
By \cite{SerreCG}*{I, \S3.3, Cor.\ 3} again, $G_2$ is torsion-free.
Furthermore, $H^*(G_1)\not\isom H^*(G_2)$.

Now if $p=2$, then $G_1$ is the free pro-$2$ product $(\dbZ/2)*\dbZ_2$, which
is an absolute Galois group of a field, e.g., an algebraic extension of $\dbQ$
\cite{Efrat99}. By Proposition \ref{principle}, $G_2$ is not a maximal pro-$p$
Galois group.

Finally, for $p>2$,  \cite{BensonLemireMinacSwallow07}*{Th.\ A.3} or
\cite{EfratMinac}*{Prop.\ 12.3} imply that in this case as well, $G_2$ is not
a maximal pro-$p$ Galois group.
\end{exam}

\begin{prop}
\label{cohomological dimension}
Let $G$ be a pro-$p$ group such that $\dim_{\dbF_p}H^1(G)<\cd(G)$.
When $p=2$ assume also that $G$ is torsion-free.
Then $G$ is not a maximal pro-$p$ Galois group of a field as above.
\end{prop}
\begin{proof}
Assume that $p\neq2$.
Let $d=\dim_{\dbF_p}H^1(G)$.
Then also $d=\dim_{\dbF_p}H^1(G^{[3,p]})$.
Since the cup product is graded-commutative, $H^{d+1}(G^{[3,p]})_\dec=0$.
On the other hand, $H^{d+1}(G)\neq0$ \cite{NeukirchSchmidtWingberg}*{Prop.\ 3.3.2}.
Thus $H^*(G^{[3,p]})_\dec\not\isom H^*(G)$.
By Theorem A and Remark \ref{www}, $G$ is not a maximal pro-$p$ Galois group.

When $p=2$ this was shown in \cite{AdemKaraMinac99}*{Th.\ 3.21}, using Kneser's theorem
on the $u$-invariant of quadratic forms, and a little later
(independently) by R.\ Ware in a letter to the third author.
\end{proof}

\begin{exam}  \rm
\label{wreath}
Let $K,L$ be finitely generated pro-$p$ groups with $1\leq n=\cd(K)<\infty$, $\cd(L)<\infty$, and $H^n(K)$ finite.
Let $\pi\colon L\to \Sym_m$, $x\mapsto \pi_x$,
be a homomorphism such that $\pi(L)$ is a transitive subgroup of $\Sym_m$.
Then $L$ acts on $K^m$ from the left by ${}^x(y_1\nek y_m)=(y_{\pi_x(1)}\nek y_{\pi_x(m)})$.
Let $G=K^m\rtimes L$.
It is generated by the generators of one copy of $K$ and of $L$.
Hence $\dim_{\dbF_p}H^1(G)=\dim_{\dbF_p}H^1(K)+\dim_{\dbF_p}H^1(L)$.

On the other hand, a routine inductive spectral sequence argument (see \cite{NeukirchSchmidtWingberg}*{Prop.\ 3.3.8})
shows that for every $i\geq0$ one has
\begin{enumerate}
\item
$\cd(K^i)=in$;
\item
$H^{in}(K^i)=H^n(K,H^{(i-1)n}(K^{i-1}))$, with the trivial $K$-action, is finite.
\end{enumerate}
Moreover, $\cd(G)=\cd(K^m)+\cd(L)$.
For $m$ sufficiently large we get $\dim_{\dbF_p}H^1(G)<mn+\cd(L)=\cd(G)$,
so by Proposition \ref{cohomological dimension}, $G$ is not a maximal pro-$p$ Galois group as above.
When $K,L$ are torsion-free, so is $G$.

For instance, one can take $K$ to be a free pro-$p$ group $\neq1$ on finitely many generators,
and let $L=\dbZ_p$ act on the direct product of $p^s$ copies of $K$ via $\dbZ_p\to\dbZ/p^s$
by cyclicly permuting the coordinates.
\end{exam}

\begin{rems}
\rm
1)\quad
In the last paragraph of the previous example, take $L=K=\dbZ_p$.
Then $p^sL\, (\isom\dbZ_p)$ acts trivially on the direct product,
so $K^{p^s}\rtimes p^sL\isom\dbZ_p^{p^s+1}$.
This group is realizable as an absolute Galois group, e.g., of $\dbC((\dbZ_{(p)}^{p^s+1}))$
\cite{EfratHaran94}*{Lemma 1.2(a)}.
Thus a torsion-free profinite group which contains an absolute Galois group of a field as an open subgroup
need not be realizable as an absolute Galois group.
This answers a question of Eli Aljadeff.

\medskip

2)\quad
Absolute Galois groups which are solvable (with respect to closed subgroups) were analyzed in
\cite{Geyer69}, \cite{Becker78}, \cite{Wurfel85}*{Cor.\ 1}, \cite{Koenigsmann01}.
In particular one can give examples of such solvable groups which are not absolute Galois groups
(compare \cite{Koenigsmann01}*{Example 4.9}).
Our examples here are in general not solvable.

\medskip

3)\quad
Another \textsl{conjectural} restriction on the structure of absolute Galois groups was suggested by Bogomolov
\cite{Bogomolov95}:
namely, let $F$ be a function field of a variety over an algebraically closed field.
Then the commutator subgroup $C$ of the $p$-Sylow subgroup of $G_F$ should be cohomologically free.
It is known that for such $F$ all elements in $H^*(C,\dbZ/p)$ are unramified (\cite{Bogomolov07}*{\S11}, 
\cite{ChernousovGilleReichstein06}*{Th.\ 1.3}).
\end{rems}

\end{document}